\documentclass[10pt]{article}
\usepackage{multicol, amsfonts, amsmath, amssymb, amsthm, setspace, epsfig, tikz, xcolor,authblk}
\usepackage[backend=biber]{biblatex}
\usepackage[margin=4cm]{geometry}
\newtheorem{thm}{Theorem}

\newtheorem{rem}{Remark}
\begin{document}
\thispagestyle{empty}

\title{Density results for the parity of \\$(4k,k)$-singular overpartitions}
\author{Victor Manuel Aricheta, Jerome Dimabayao, Hazel Joy Shi}
\date{}

\maketitle

\begin{abstract}
The $(k,i)$-singular overpartitions, combinatorial objects introduced by Andrews in 2015, are known to satisfy Ramanujan-type congruences modulo any power of prime coprime to $6k$. In this paper we consider the parity of the number $\overline{C}_{k,i}(n)$ of $(k,i)$-singular overpartitions of $n$. In particular, we give a sufficient condition on even values of $k$ so that the values of $\overline{C}_{4k,k}(n)$ are almost always even. Furthermore, we show that for odd values of $k \leq 23$, $k\neq 19$, certain subsequences of $\overline{C}_{4k,k}(n)$ are almost always even.
\end{abstract}

\section{Introduction}

Singular overpartitions were introduced by Andrews in \cite{andrews2015singular}, with the primary motivation of presenting a theorem for overpartitions that is analogous to Gordon's combinatorial generalization of the Rogers-Ramanujan identities \cite{gordon1961combinatorial}. In the same paper and in addition to this combinatorial property, Andrews showed that singular overpartitions satisfy congruences reminiscent of Ramanujan's celebrated congruences for the ordinary partition function \cite{ramanujan1921congruence}. It is now known that, as in the case of the ordinary partition function \cite{ono2000distribution}, there are infinitely many Ramanujan-type congruences for singular overpartitions for all but finitely many prime moduli \cite{aricheta2017congruences}. The exceptional primes not covered by this result include the prime 2. In this paper we establish some density results on the parity of singular overpartitions.

\subsection{Partitions and singular overpartitions}

A \emph{partition} of a positive integer $n$ is a non-increasing sequence of positive integers that sum up to $n$. The number of partitions of a positive integer $n$ is denoted by $p(n)$, and we set $p(0) = 1$. Partitions have been studied as far back by Euler, who obtained the following generating function
\[ \sum_{n=0}^\infty p(n) q^n = \prod_{m=1}^\infty \dfrac{1}{1-q^m}.\]
Ramanujan showed that the partition function $p(n)$ satisfies congruences of the form $p(An+B) \equiv 0 \pmod \ell$ for $\ell = 5,7,11$. Using modular forms, Ahlgren and Ono showed that there are many such congruences \cite{ ahlgren2000distribution, ono2000distribution}, that is, for each positive integer $\ell$ coprime to 6, there are infinitely many non-nested arithmetic progressions $An+B$ for which $p(An+B) \equiv 0 \pmod \ell$.

Less is known about the behavior of $p(n)$ modulo $2$ or modulo $3$. We know from the works of Ono and Radu that given positive integers $r$ and $t$, there are infinitely many $N \equiv r \pmod t$ such that $p(N)$ is even \cite{ono1996parity}, and infinitely many $M \equiv r \pmod t$ such that $p(M)$ is odd \cite{radu2012proof}. This supports the conjecture \cite{parkin1967distribution} that the sequence $p(n)$ is equidistributed modulo 2 and modulo 3, i.e., that:
\[ \lim_{M\to \infty} \dfrac{\# \{1 \leq n\leq M : p(n) \equiv 0 \text{ (mod  2)}\}}{M} = \dfrac{1}{2}, \]
\[ \lim_{M\to \infty} \dfrac{\# \{1 \leq n\leq M : p(n) \equiv 0 \text{ (mod  3)}\}}{M} = \dfrac{1}{3}. \]

An \emph{overpartition} of $n$ is a non-increasing sequence of positive integers whose sum is $n$ and in which the first occurrence of a number may be overlined. Overpartitions are generalizations of partitions, introduced (and developed in a series of papers) by Lovejoy and Corteel \cite{corteel2004overpartitions}. Andrews showed that overpartitions satisfy a theorem analogous to Gordon's generalization of the Rogers-Ramanujan identities \cite{andrews2015singular}. To do so, he introduced the notion of \emph{$(k,i)$-singular overpartitions}. We dispense with the exact definition of singular overpartitions due to its involved nature. (We refer the reader to \cite{andrews2015singular} for the exact definition.) The following result, also proven by Andrews, is sufficient for our needs: Let $\overline{C}_{k,i}(n)$ be the number of $(k,i)$-singular overpartitions of $n$. Then $\overline{C}_{k,i}(n)$ is equal to the number of overpartitions of $n$ in which no part is divisible by $k$ and only the parts that are congruent to $\pm i \pmod k$ may be overlined.  We set $\overline{C}_{k,i}(0) = 1$. We then have the following generating function for $\overline{C}_{k,i}(n)$:
\[ \sum_{n=0}^\infty \overline{C}_{k,i}(n) q^n = \prod_{m = 1}^\infty \frac{\left(1 - q^{km}\right)\left(1 + q^{k(m - 1) + i}\right)\left(1 + q^{km - i}\right)}{1 - q^m}.\]

\subsection{Statement of results}

In this paper we consider the parity of $(k,i)$-singular overpartitions. This problem has been considered in many recent papers \cite{aricheta2017congruences,barman2019divisibility,liu2019arithmetic,singh2021new}. In particular, the first author proved that $ \overline{C}_{3k, k}(n) $ is ``almost always even'' whenever $k = 2^am$ where $a$ is a nonnegative integer, $m$ is a positive odd number, and $ 2^a \geq m $; and Singh and Barman proved that $ \overline{C}_{4k, k} (n)$ is ``almost always divisible by arbitrary powers of $2$'' with analogous conditions for $k$ as in \cite{aricheta2017congruences}. We note here that except for $k = 1$, Singh and Barman's result covers only even values of $k$.

In particular we consider $(4k,k)$-singular overpartitions and show that density results modulo 2 are already quite varied. These $(4k,k)$-singular overpartitions modulo 2 are known to be related to the minimal excludant of partitions \cite{barman2020mex}. 

Before we state our results, we first make precise the notion of being ``almost always even.'' Consider a positive integer $m$. A power series $ \sum_{n = 0}^\infty a_n q^n $ is said to be  \emph{lacunary modulo $m$} if
\[ \lim_{x \to \infty} \dfrac{\#\{a_n : n \leq x, a_n \equiv 0 \pmod{m}\}}{x}  = 1.\]
If this is the case, we also say that \emph{the sequence $a_n$ is lacunary modulo $m$}. In particular, if the sequence $a_n$ is lacunary modulo $2$, then the values of $a_n$ are ``almost always even.''

We first consider even values of $k$. Write $k = 2^a\ell$ where $a, \ell$ are positive integers and $\ell$ is odd. In \cite{singh2021new}, Singh and Barman showed that if $\ell \leq 2^a$, then the sequence $\overline{C}_{4k,k}(n)$ is lacunary modulo arbitrary powers of 2. The following theorem---which we obtain by an application of the results of Cotron et al.\ \cite{cotron2020lacunary} on eta-quotients---extends the range of $\ell$ for which the same conclusion holds, at the expense of restricting the modulus to 2.

\begin{thm}
Let $k$ be a positive even integer. Write $k = 2^a \ell$ where $a,\ell$ are positive integers and $\ell$ is odd. If $\ell \leq 3\cdot 2^a$ then the sequence $\overline{C}_{4k,k}(n)$ is lacunary modulo 2.
\end{thm}

\begin{rem} The condition $\ell \leq 3\cdot 2^a$ in Theorem 1 appears to be sharp; numerical computations suggest that the sequence $\overline{C}_{4k,k}(n)$ is equidistributed modulo 2 if $\ell > 3\cdot 2^a$. To illustrate, let
\[ \delta_k(M) := \dfrac{\# \{1 \leq n\leq M : \overline{C}_{4k,k}(n) \equiv 0 \text{ (mod  2)}\}}{M} \]
be the proportion of even $(4k,k)$-singular overpartition numbers among the first $M$ positive integers. The following table gives some values of $\delta_k(M)$ for $k = 2\ell$ for odd $\ell \leq 11$. 
\begin{center}
\begin{tabular}{|c|c|c|c|c|c|c|}
\hline
$M$ & $\delta_{2}(M)$ & $\delta_{6}(M)$ & $\delta_{10}(M)$ & $\delta_{14}(M)$ & $\delta_{18}(M)$ & $\delta_{22}(M)$\\
\hline
$100$ & $0.46$ & 0.5 & 0.51 & 0.47 & 0.47 & 0.39 \\
$1000$ & $0.603$ & 0.607 & 0.617 & 0.48 & 0.477 & 0.477\\
$10000$ & $0.688$ & 0.6923 & 0.6927 & 0.4937 & 0.5014 & 0.4986 \\
$100000$ & $0.7455$ & 0.74656 & 0.74705 & 0.499 & 0.50012 & 0.49813 \\
$\vdots$ & $\vdots$ & $\vdots$ & $\vdots$ & $\vdots$ & $\vdots$ & $\vdots$\\
$\infty$ & 1 & 1 & 1 & 0.5? & 0.5? & 0.5? \\
\hline
\end{tabular}
\end{center}
Note that for odd integers $\ell \leq 6$, Theorem 1 guarantees that $\delta_{2\ell}(M)$ approaches 1 as $M$ tends to infinity. On the other hand, for odd integers $\ell > 6$, the values of $\overline{C}_{4\cdot 2\ell,2\ell}(n)$ appear to be equidistributed modulo 2.
As another example, consider $k = 4\ell$ with $\ell$ odd. Then by Theorem 1, the sequence $\overline{C}_{4k,k}(n)$ is lacunary modulo 2 when $\ell \leq 12$. In contrast, the values of $\delta_{4\ell}(100000)$ for $\ell=13,15,17$ are 0.49684, 0.49725 and 0.49929, respectively, and these too are close to 1/2. 
\end{rem}

We now turn our attention to the case when $k$ is odd, and in particular look at the cases $k\leq 23$, $k \neq 19$. We will prove the following.

\begin{thm} The following sequences are lacunary modulo 2:
\begin{enumerate} 
\item ($k=1$) the sequence $\overline{C}_{4,1}(n)$,
\item ($k=3$) the sequence $\overline{C}_{12,3}(n)$,
\item ($k=5$) the sequence $\overline{C}_{20,5}(2n)$,
\item ($k=7$) the sequences $\overline{C}_{28,7}(2n+1)$ and $\overline{C}_{28,7}(4n)$,
\item ($k=9$) the sequence $\overline{C}_{36,9}(2n+1)$,
\item ($k=11$) the sequences $\overline{C}_{44,11}(2n)$ and $\overline{C}_{44,11}(4n+1)$,
\item ($k=13$) the sequence $\overline{C}_{52,13}(2n)$,
\item ($k=15$) the sequences $\overline{C}_{60,15}(2n+1)$ and $\overline{C}_{60,15}(4n+2)$,
\item ($k=17$) the sequence $\overline{C}_{68,17}(2n+1)$,
\item ($k=21$) the sequence $\overline{C}_{84,21}(2n)$,
\item ($k=23$) and the sequence $\overline{C}_{92,23}(2n+1)$.
\end{enumerate}
\end{thm}

From these small values of odd $k$, we already find that the results are more varied and the situation more complex than the case when $k$ is even. Numerical evidence indicates that the following sequences are equidistributed modulo 2: $\overline{C}_{20,5}(2n+1)$, $\overline{C}_{28,7}(4n+2)$, $\overline{C}_{36,9}(2n)$, $\overline{C}_{52,13}(2n+1)$, $\overline{C}_{60,15}(4n)$, $\overline{C}_{68,17}(2n)$ and $\overline{C}_{84,21}(2n+1)$. These are precisely the sequences excluded by Theorem 2 for $k=5,7,9,13,17$ and $21$. Thus, in conjunction with the results of Theorem 2, we conjecture that $\overline{C}_{4k,k}(n)$ is even $3/4$ of the time for $k = 5,9,13,17$ and $21$, while $\overline{C}_{4k,k}(n)$ is even $7/8$ of the time for $k=7$ and $15$. We note that the generating function for $\overline{C}_{4k,k}(n)$ is an eta-quotient (cf. Section 2.1), and our parity conjectures are consistent with the general conjecture on the density of odd (or even) coefficients of eta-quotients recently proposed by Keith and Zanello \cite{keith2022parity}.

\vspace{0.1in}

\noindent \textbf{Acknowledgment}. The authors thank Richell Celeste for her helpful suggestions and comments regarding Theorem 1, and the anonymous referee for their valuable suggestions which greatly expanded the cases covered by Theorem 2.

\section{Preliminaries}

In this paper we use the notation
\[ f_k := \prod_{m = 1}^\infty (1-q^{km}).\]
Since all the congruences in this paper are congruences modulo 2, we simply write 
\[ A \equiv B\]
as a shorthand for $A \equiv B \pmod 2$.

Set $ \overline{C}_{k, i}(0) = 1$. The generating function for $ \overline{C}_{k, i}(n) $ is
\[ \sum_{n = 0}^\infty \overline{C}_{k, i}(n)q^n = \prod_{m = 1}^\infty \frac{\left(1 - q^{km}\right)\left(1 + q^{k(m - 1) + i}\right)\left(1 + q^{km - i}\right)}{1 - q^m}. \]
In particular the generating function for $ \overline{C}_{4k, k}(n) $ is
\begin{align} & \sum_{n = 0}^\infty \overline{C}_{4k, k}(n)q^n \nonumber \\
	& = \prod_{m = 1}^\infty \frac{\left(1 - q^{4km}\right)\left(1 + q^{4km - 3k}\right)\left(1 + q^{4km - k}\right)}{1 - q^m} \cdot \frac{\left(1 + q^{4km - 2k}\right)\left(1 + q^{4km}\right)}{\left(1 + q^{4km - 2k}\right)\left(1 + q^{4km}\right)} \nonumber \\
	& = \prod_{m = 1}^\infty \frac{\left(1 - q^{4km}\right)\left(1 + q^{km}\right)}{\left(1 - q^m\right)\left(1 + q^{2km}\right)} \cdot \frac{\left(1 + q^{km}\right)\left(1 - q^{km}\right)^2}{\left(1 + q^{km}\right)\left(1 - q^{km}\right)^2} \nonumber \\
	& = \prod_{m = 1}^\infty \frac{\left(1 - q^{4km}\right)\left(1 - q^{2km}\right)^2}{\left(1 - q^m\right)\left(1 - q^{km}\right)\left(1 - q^{4km}\right)} \nonumber \\
	& \equiv \prod_{m = 1}^\infty \frac{\left(1 - q^{km}\right)^4}{\left(1 - q^m\right)\left(1 - q^{km}\right)}\nonumber \\
	& \equiv \prod_{m = 1}^\infty \frac{\left(1 - q^{km}\right)^3}{\left(1 - q^m\right)} \nonumber \\
	&  \equiv \dfrac{f_k^3}{f_1}. \label{genfcn}
\end{align}
Let $\tau \in \mathbb{H} := \{ \tau \in \mathbb{C} : \text{Im}(\tau) > 0 \}$ and let $q = e^{2\pi i \tau} $. The Dedekind eta-function is the function
 \[ \eta(\tau) := q^\frac{1}{24} \prod_{m = 1}^\infty \left(1 - q^m\right).\]
Then the generating function for $ \overline{C}_{4k, k}(n)$ is also given by
\begin{equation}\label{eta}
\sum_{n = 0}^\infty \overline{C}_{4k, k}(n)q^n \equiv q^{\frac{1 - 3k}{24}} \dfrac{\eta(k\tau)^3}{\eta(\tau)}.
\end{equation}

We gather here several theorems that we will use in proving that certain functions are lacunary modulo 2. Theorem \ref{Cotron} is due to Cotron et al.\ \cite{cotron2020lacunary}, while Theorem \ref{Landau} is due to Landau \cite{landau1909einteilung}.

\begin{thm} \label{Cotron}
Let $s$ and $t$ be positive integers with $s$ odd and $s < t$. Let $ p $ be a prime divisor of $ t $ and let $ a $ be a positive integer. If $ p^a \mid t $ and  $p^a \geq \sqrt{\frac{t}{s}}$,
then $ E_{t, s}(z) := q^{\frac{1 - ts}{24}} \frac{\eta^s (tz)}{\eta(z)} $ is lacunary modulo $ p^j $ for any positive integer $ j $.
\end{thm}

\begin{thm}[Landau's Theorem] \label{Landau}
Let $ r(m) $ and $ s(n) $ be quadratic polynomials on $ m $ and $ n $, respectively, then
\[ \left(\sum_{m \in \mathbb{Z}} q^{r(m)}\right)\left(\sum_{n \in \mathbb{Z}} q^{s(n)}\right) \]
is lacunary modulo $2$.
\end{thm}

We also recall Euler's pentagonal number theorem modulo 2:
\begin{equation}\label{euler}
f_1 \equiv \sum_{r \in \mathbb{Z}} q^{r(3r - 1)/2},
\end{equation}
from which it follows for a positive integer $t$ that
\begin{equation}\label{euler-t}
f_t \equiv \sum_{r \in \mathbb{Z}} q^{tr(3r - 1)/2}.
\end{equation}
Finally we have a special case of the Jacobi triple product modulo 2:
\begin{equation}\label{jacobi}
f_1^3 \equiv \sum_{s \in \mathbb{Z}} q^{s(s + 1)/2},
\end{equation}
from which we obtain that
\begin{equation}\label{jacobi-t}
f_t^3 \equiv \sum_{s \in \mathbb{Z}} q^{ts(s + 1)/2},
\end{equation}
for a positive integer $t$.
\section{Proof of Theorem 1}

\begin{proof}
We first consider the case when $k=2$, i.e., the case when $a = \ell = 1$. From congruence (\ref{genfcn}) we have
\begin{align*} \sum_{n = 0}^\infty \overline{C}_{8, 2}(n) q^n & \equiv \dfrac{f_2^3}{f_1}\\
		& \equiv f_1^5\\
		& \equiv f_1^2 f_1^3.
\end{align*}
Using congruences (\ref{euler}) and (\ref{jacobi}) we get
\[ \sum_{n = 0}^\infty \overline{C}_{8, 2}(n) q^n \equiv \left(\sum_{r \in \mathbb{Z}} q^{r(3r - 1)}\right)\left(\sum_{s \in \mathbb{Z}} q^{(s^2 + s)/2}\right). \]
We conclude by Landau's theorem that the sequence $\overline{C}_{8, 2}(n) $ is lacunary modulo $ 2$.

We now consider the case when $k = 2^a \ell$, where $ a, \ell$ are positive integers satisfying the following: $\ell$ is odd, the pair $(a,\ell) \neq (1,1)$, and $\ell \leq 3 \cdot 2^a$. From congruence (\ref{eta}), we have
\begin{align*}
\sum_{n = 0}^\infty \overline{C}_{4k, k}(n)q^n &\equiv q^{\frac{1 - 3k}{24}} \frac{\eta^3 (k\tau)}{\eta(\tau)} \\
& \equiv E_{k, 3}(\tau)
\end{align*}
where $E_{t, s}(\tau)$ is the eta-quotient appearing in Theorem \ref{Cotron}. Note that $ k \geq 4 > 3 $ and $ 2^a \mid k $. Furthermore, since $\ell \leq 3 \cdot 2^a$, we have
\[ \sqrt{\frac{k}{3}} = \sqrt{\frac{2^a\ell}{3}} \leq \sqrt{\frac{3 \cdot 2^{2a}}{3}} = 2^a. \]
Therefore by Theorem \ref{Cotron} (with $t = k$, $s = 3$ and $p = 2$) the eta-quotient $E_{k, 3}(\tau)$ and the sequence $\overline{C}_{4k, k}(n)$ are lacunary modulo $ 2 $.
\end{proof}

\section{Proof of Theorem 2}

\subsection{Case $k = 1$}
\begin{proof}
We know from \cite[Theorem 3.2]{chen2015arithmetic} that $\overline{C}_{4, 1}(n)$ is odd if and only if $n = m(3m-1)$ for some integer $m$. Thus the sequence $ \overline{C}_{4, 1}(n) $ is lacunary modulo $ 2 $.
\end{proof}

\subsection{Case $k = 3$}
\begin{proof}
From congruence (\ref{genfcn}) we have
\[ \sum_{n = 0}^\infty \overline{C}_{12, 3}(n) q^n \equiv \dfrac{f_3^3}{f_1}.\]
The right hand side of this congruence is the generating function for the number $c_3(n)$ of $3$-core partitions of $n$, and we have the congruence (cf.\ \cite[p.\ 266]{ono1997odd})
\begin{equation}\label{3core}
\dfrac{f_3^3}{f_1} \equiv \sum_{m\in \mathbb{Z}} q^{(3m^2 + 2m)/2}.
\end{equation}
Therefore the sequence $\overline{C}_{12, 3}(n)$ is lacunary modulo 2.
\end{proof}

\subsection{Case $k=5$}

\begin{proof} Begin with the even-odd dissection of $\frac{f_5}{f_1}$ (cf. \cite[Section 8]{keith2022parity})
\begin{equation}\label{f5f1}
\dfrac{f_5}{f_1} \equiv f_1^4 + q\dfrac{f_5^6}{f_1^2}
\end{equation}
and multiply by $f_5^2$ to obtain
\[ \sum_{n = 0}^\infty \overline{C}_{20, 5}(n) q^{n} \equiv \dfrac{f_5^3}{f_1} \equiv f_1^4 f_5^2 + q\dfrac{f_5^8}{f_1^2}. \] 
By extracting the terms with even powers of $q$ we get
\begin{align*}
\sum_{n = 0}^\infty \overline{C}_{20, 5}(2n) q^{2n} &\equiv f_1^4f_5^2\\
&\equiv f_2^2f_5^2.
\end{align*}
Replacing $q$ by $q^{1/2}$ yields
\[ \sum_{n = 0}^\infty \overline{C}_{20, 5}(2n) q^{n} \equiv f_2 f_5. \]
Using congruence (\ref{euler-t}) we obtain
\[ \sum_{n = 0}^\infty \overline{C}_{20, 5}(2n) q^{n} \equiv \left(\sum_{r \in \mathbb{Z}} q^{r(3r - 1)}\right)\left(\sum_{r \in \mathbb{Z}} q^{\frac{5r(3r - 1)}{2}}\right). \]
Therefore, by Landau's theorem, the sequence $ \overline{C}_{20, 5}(2n)$ is lacunary modulo $2$.
\end{proof}

\subsection{Case $k=7$}

\subsubsection{Lacunarity of $\overline{C}_{28, 7}(2n+1)$}

\begin{proof} We start with the following congruence given in \cite[Congruence (2.4)]{lin2014elementary}:
\[ f_1f_7 \equiv f_1^8 + qf_1^4f_7^4 + q^2 f_7^8 \]
and multiply this by $f_7^2/f_1^2$ to obtain
\[ \sum_{n = 0}^\infty \overline{C}_{28, 7}(n) q^{n} \equiv \dfrac{f_7^3}{f_1} \equiv f_1^6 f_7^2 + qf_1^2f_7^6 + q^2\dfrac{f_7^{10}}{f_1^2}\] 
Extracting the terms with odd powers of $q$, then dividing both sides by $q$ and then replacing $q$ by $q^{1/2}$ yields
\[ \sum_{n = 0}^\infty \overline{C}_{28, 7}(2n+1) q^{n} \equiv f_1f_7^3. \]
By congruences (\ref{euler}) and (\ref{jacobi-t}) and Landau's theorem, the sequence $\overline{C}_{28, 7}(2n + 1)$ is lacunary modulo $2$.
\end{proof}

\subsubsection{Lacunarity of $\overline{C}_{28, 7}(4n)$}

\begin{proof}
Let $c_7(n)$ be the number of 7-core partitions of $n$. Then
\begin{align*}
 \sum_{n = 0}^\infty c_7(n) q^{n} &\equiv \dfrac{f_7^7}{f_1}\\
 & \equiv \dfrac{f_7^3}{f_1}f_{28}\\
 & \equiv \left(\sum_{n = 0}^\infty \overline{C}_{28, 7}(n) q^{n}\right)f_{28}
\end{align*}
Consider the operator $U_2$ defined as follows:
\[ U_2\left(\sum_{n=0}^\infty a(n)q^n\right) = \sum_{n=0}^\infty a(2n)q^n. \]
It is straightforward to check, by expanding both sides of the equation, that
\[ U_2\left(\sum_{n = 0}^\infty a(n) q^n \sum_{m = 0}^\infty b(m) q^{2m}\right) = \left(\sum_{n = 0}^\infty a(2n) q^n\right)\left(\sum_{m = 0}^\infty b(m) q^m\right).\]
Therefore, applying $U_2$ to the generating function of $c_7(n)$ twice we get:
\[
\sum_{n = 0}^\infty c_7(4n) q^{n} \equiv \left(\sum_{n = 0}^\infty \overline{C}_{28, 7}(4n) q^{n}\right) f_{7}.
\]
From \cite[Theorem 2.9]{xia2014new} we have that
\[ \sum_{n = 0}^\infty c_7(4n) q^n \equiv f_1 f_4 f_7, \]
and therefore, from the previous two congruences, we get
\[ \sum_{n = 0}^\infty \overline{C}_{28, 7}(4n) q^n \equiv f_1 f_4 \]
By congruences (\ref{euler}) and (\ref{euler-t}) and Landau's theorem, we get that $ \overline{C}_{28, 7}(4n) $ is lacunary modulo $2$.
\end{proof}

\subsection{Case $k=9$}

\begin{proof} Begin with the following even-odd dissection of $f_9/f_1$ (cf.\ \cite[Equation (2.1)]{yao2014new}:
\[ \dfrac{f_9}{f_1} = \dfrac{f_{12}^3 f_{18}}{f_2^2 f_6 f_{36}} + q \dfrac{f_4^2 f_6 f_{36}}{f_2^3 f_{12}} \] 
and multiply this by $f_9^2$ to get
\[ \sum_{n = 0}^\infty \overline{C}_{36, 9}(n) q^{n} \equiv \dfrac{f_{9}^3}{f_1} = \dfrac{f_{12}^3 f_{18}f_9^2}{f_2^2 f_6 f_{36}} + q \dfrac{f_4^2 f_6 f_{36}f_9^2}{f_2^3 f_{12}}. \] 
Extracting the terms with odd powers of $q$ we obtain
\begin{align*} \sum_{n = 0}^\infty \overline{C}_{36, 9}(2n+1) q^{2n+1} &\equiv q \dfrac{f_4^2 f_6 f_{36} f_9^2}{f_2^3 f_{12}}\\
&\equiv q\dfrac{f_1^2 f_9^6}{f_3^2}.
\end{align*}
Dividing by $q$ and replacing $q$ by $q^{1/2}$ yields
\begin{align*}
\sum_{n = 0}^\infty \overline{C}_{36, 9}(2n+1) q^{n} &\equiv f_1 \dfrac{f_9^3}{f_3} \\
&\equiv \left(\sum_{r \in \mathbb{Z}} q^{r(3r - 1)/2}\right)\left( \sum_{m\in \mathbb{Z}} q^{3(3m^2 + 2m)/2}\right)
\end{align*}
where the congruence for $f_9^3/f_3$ is obtained from replacing $q$ by $q^3$ in congruence (\ref{3core}).
The lacunarity of $\overline{C}_{36, 9}(2n+1)$ follows from using Landau's theorem.
\end{proof}

\subsection{Case $k=11$}

Begin with the following congruence (cf. \cite[Lemma 2.1]{zhao2018parity}): 
\[
\frac{1}{f_1  f_{11}} 
\equiv 
\frac{f_{12}^3}{f_4 f_{44}} + 
q \frac{f_6^3}{f_2^3} + q^6 \frac{f_{66}^3}{f_{22}^3} + q^{15} \frac{f_{132}^3}{f_4 f_{44}}
\]
and multiply this by $f_{11}^4$ to get 
\begin{equation}\label{k is 11}
\sum_{n = 0}^\infty \overline{C}_{44, 11}(n) q^n \equiv 
\frac{f_{11}^3}{f_1} 
\equiv 
\frac{f_{12}^3}{f_4} + 
q \frac{f_6^3 f_{11}^4}{f_2^3} + q^6 \frac{f_{66}^3}{f_{22}} + q^{15} \frac{f_{132}^3}{f_4}.
\end{equation}

\subsubsection{Lacunarity of $\overline{C}_{44, 11}(2n)$}

\begin{proof}
Extracting the terms in (\ref{k is 11}) with even powers of $q$ and replacing $q$ by $q^{1/2}$ yields 
\[
\sum_{n = 0}^\infty \overline{C}_{44, 11}(2n) q^{n} \equiv 
\frac{f_{6}^3}{f_2} + q^3 \frac{f_{33}^3}{f_{11}} \pmod{2}.
\]
Since $f_{6}^3/f_2$ and $f_{33}^3/f_{11}$ are just the $q \rightarrow q^2$ and $q \rightarrow q^{11}$ magnifications of 
$f_3^3/f_1$ respectively, it follows from congruence (\ref{3core}) that each of the summands above are lacunary modulo $2$. Therefore, the sequence $\overline{C}_{44, 11}(2n)$ is lacunary modulo $2$.
\end{proof}

\subsubsection{Lacunarity of $\overline{C}_{44, 11}(4n+1)$}

\begin{proof}
We extract the terms in (\ref{k is 11}) with odd powers of $q$, then divide by $q$, and then
replace $q$ by $q^{1/2}$ to obtain 
\begin{equation}\label{odd-11}
 \sum_{n = 0}^\infty \overline{C}_{44, 11}(2n+1) q^{n}
\equiv 
\left( \frac{f_3}{f_1} \right)^{3}  f_{11}^2  + q^{7} \frac{f_{66}^3}{f_2}. 
\end{equation}
From the following even-odd dissection of $f_3/f_1$ (cf.\ \cite[Congruence (6)]{keith2022parity}): 
\begin{equation}\label{3regular}
\frac{f_3}{f_1} \equiv \frac{f_1^8}{f_3^2} + q \frac{f_3^{10}}{f_1^4}
\end{equation}
we have
\begin{equation*}
\left( \frac{f_3}{f_1} \right)^{3} \equiv 
\left( \frac{f_1^8}{f_3^2}  + q \frac{f_3^{10}}{f_1^{4}}\right)^3  
\equiv 
\frac{f_1^{24}}{f_3^6} + q f_1^{12} f_3^{6}  + q^2 f_{3}^{18} + q^3 \frac{f_3^{30}}{f_1^{12}}.
\end{equation*}
We substitute this into relation (\ref{odd-11}), then extract the terms with even powers of $q$, 
and then replace $q$ by $q^{1/2}$ to get 
\[
 \sum_{n \geq 0} \overline{C}_{44, 11}(4n+1) q^{n} 
 \equiv \left( \frac{f_1^{12}}{f_3^{3}} + q f_3^{9} \right) f_{11} 
 \equiv 
 \left( \frac{f_1^{12} + q f_3^{12}} {f_3^{3}} \right) f_{11}.
\]
We then use the congruence (cf.\ \cite[Congruence (5)]{keith2022parity})
\begin{equation}\label{f1f3}
(f_1f_3)^3 \equiv f_1^{12} + qf_3^{12}
 \end{equation}
in place of the numerator inside the parenthesis to get
\[
 \sum_{n = 0}^\infty \overline{C}_{44, 11}(4n+1) q^{n} 
 \equiv f_1^3 f_{11}.
\]
By congruences (\ref{euler-t}) and (\ref{jacobi}) and Landau's theorem, the sequence $\overline{C}_{44, 11}(4n+1)$ is lacunary modulo $2$.
\end{proof}

\subsection{Case $k=13$}

\begin{proof} Begin with the following congruence  (cf.\ \cite[p. 589]{judge2018density}):
\[ \dfrac{f_{13}}{f_1} \equiv f_1^{12} + qf_1^{10} f_{13}^2 + q^6f_{13}^{12} + q^7\dfrac{f_{13}^{14}}{f_1^2} \] 
and multiply this by $f_{13}^2$ to get
\[ \sum_{n = 0}^\infty \overline{C}_{52, 13}(n)q^n \equiv \dfrac{f_{13}^3}{f_1} \equiv f_1^{12}f_{13}^2 + qf_1^{10} f_{13}^4 + q^6f_{13}^{14} + q^7\dfrac{f_{13}^{16}}{f_1^2}. \] 
Extracting the terms with even powers of $q$, and replacing $q$ by $q^{1/2}$ yields
\begin{align*}
\sum_{n = 0}^\infty \overline{C}_{52, 13}(2n) q^{n} & \equiv f_1^{6}f_{13} + q^3f_{13}^{7}\\
		& \equiv f_2^3 f_{13} + q^3f_{26}^3 f_{13}
\end{align*}
Using congruences (\ref{euler-t}) and (\ref{jacobi-t}), and Landau's theorem, we find that each of the summand on the right hand side is lacunary modulo 2. The sum of two lacunary series is also lacunary, therefore the sequence $ \overline{C}_{52, 13}(2n) $ is also lacunary modulo $2$.
\end{proof}

\subsection{Case $k=15$}

Divide the congruence (\ref{f1f3}) by $f_1^4$, and replace $q$ by $q^5$ to get
\begin{equation}\label{Rel for k=3 magnified by 5}
\frac{f_{15}^3}{f_5} \equiv f_5^8 + q^5 \frac{f_{15}^{12}}{f_5^4}.
\end{equation}
From congruences (\ref{f5f1}) and (\ref{Rel for k=3 magnified by 5}) we find that 
\begin{align*}
\sum_{n = 0}^\infty \overline{C}_{60, 15}(n)q^n  & \equiv \frac{f_{15}^3}{f_5} \frac{f_5}{f_1} \\
& \equiv 
\left( f_5^8 + q^5 \frac{f_{15}^{12}}{f_5^4} \right)\left( f_1^4 + q \frac{f_5^6}{f_1^2} \right) \\
&\equiv 
f_1^4 f_5^8 + q \frac{f_5^{14}}{f_1^2} + q^5 \frac{f_1^4 f_{15}^{12}}{f_5^4} + q^6 \frac{f_5^2 f_{15}^{12}}{f_1^{2}}.
\end{align*}
Thus,
\begin{equation}\label{k=15,even}
\sum_{n = 0}^\infty \overline{C}_{60, 15}(2n)q^n \equiv f_1^2 f_5^4 + q^3 \frac{f_5 f_{15}^6}{f_1}
\end{equation}
and 
\begin{equation}\label{k=15,odd}
\sum_{n = 0}^\infty \overline{C}_{60, 15}(2n+1)q^n \equiv \frac{f_5^{7}}{f_1} + q^2 \frac{f_1^2 f_{15}^{6}}{f_5^2}.
\end{equation}

\subsubsection{Lacunarity of $\overline{C}_{60, 15}(2n+1)$}
\begin{proof}
We write relation (\ref{k=15,odd}) as 
\begin{align*}
\sum_{n = 0}^\infty \overline{C}_{60, 15}(2n+1)q^n &\equiv
\dfrac{f_5^5}{f_1} f_{10} + q^2 f_2 \dfrac{f_{30}^3}{f_{10}}.
\end{align*}
The quotient $f_{30}^3/f_{10}$ is the $q \rightarrow q^{10}$ magnification of $f_3^3/f_1$, and therefore by congruences (\ref{euler-t}) and (\ref{3core}) and Landau's theorem, the second term in the above sum is lacunary modulo $2$. To deal with the first term in the above sum, we use the following congruence (cf.\ \cite[Lemma 18]{keith2022parity}): 
\[ \dfrac{f_5^5}{f_1} \equiv \sum_{n = 1}^\infty q^{n^2-1} + \sum_{n = 1}^\infty q^{2n^2-1} + \sum_{n = 1}^\infty q^{5n^2-1} +  \sum_{n = 1}^\infty q^{10n^2-1}.\]
From this and congruence (\ref{euler-t}) we see that
\[ \dfrac{f_5^5}{f_1} f_{10} \equiv \left(\sum_{n = 1}^\infty q^{n^2-1}  + \sum_{n = 1}^\infty q^{2n^2-1} + \sum_{n = 1}^\infty q^{5n^2-1} +  \sum_{n = 1}^\infty q^{10n^2-1}\right) \left( \sum_{r \in \mathbb{Z}} q^{5r(3r - 1)} \right).\]
By Landau's theorem, each of the four terms here is lacunary modulo $2$. Therefore, the sequence $\overline{C}_{60, 15}(2n+1)$ is also lacunary modulo $2$.
\end{proof}

\subsubsection{Lacunarity of $\overline{C}_{60, 15}(4n+2)$}
\begin{proof}
By congruence (\ref{f5f1}) we may write congruence
(\ref{k=15,even}) as
\[
\sum_{n = 0}^\infty \overline{C}_{60, 15}(2n)q^n 
\equiv f_1^2 f_5^4 + q^3 f_{15}^6 \left(f_1^4 + q \frac{f_5^6}
{f_1^2} \right).
\]
We extract the terms with odd powers of $q$, then divide by $q$ and then replace
$q$ by $q^{1/2}$ to get
\[
\sum_{n = 0}^\infty \overline{C}_{60, 15}(4n+2)q^{n} \equiv q f_{15}^3
f_2.
\]
By congruences (\ref{euler-t}) and (\ref{jacobi-t}) and Landau's theorem, the sequence $\overline{C}_{60, 15}(4n+2)$ is lacunary modulo $2$.
\end{proof}

\subsection{Case $k=17$}

\begin{proof}
Begin by considering congruence (4.5) in \cite{zhao2018parity}: 
\[ \frac{1}{f_1 f_{17}} \equiv \sum_{n = 0}^\infty \Delta_8(2n) q^{2n} + q f_2^3 + q^5 f_{34}^4,\]
where $\Delta_8(n)$ is defined to be the coefficient of $q^n$ in $f_2 f_{17}/f_1^3 f_{34}$. Multiplying this by $f_{17}^4$ yields
\[ \sum_{n = 0}^\infty \overline{C}_{68,17}(n) q^n \equiv 
\frac{f_{17}^3}{f_1} \equiv f_{17}^4 \sum_{n \geq 0} \Delta_8(2n) q^{2n} + q f_2^3 f_{17}^4  + q^5 f_{17}^{12}. \]
Extract the terms with odd powers of $q$, then divide by $q$, then replace $q$ by $q^{1/2}$ to obtain 
\[ \sum_{n \geq 0} \overline{C}_{68,17}(2n+1) q^{n} \equiv f_1^3 f_{34}  + q^2 f_{34}^{3}.\]
The term $f_1^3 f_{34}$ is lacunary modulo $2$ by congruences (\ref{euler-t})  and (\ref{jacobi}) and Landau's theorem, while the term $f_{34}^{3}$ is lacunary modulo $2$ by congruence (\ref{jacobi-t}). Hence the sequence $\overline{C}_{68,17}(2n+1)$ is lacunary modulo $2$.
\end{proof}

\subsection{Case $k=21$}

\begin{proof}
Replacing $q$ by $q^7$ in congruence (\ref{3regular}), and multiplying both sides of the resulting congruence by $f_{21}^2$ gives us
\begin{equation}\label{relfork=21(1)}
\frac{f_{21}^3}{f_7} \equiv f_7^8 + q^7 \frac{f_{21}^{12}}{f_7^4}. 
\end{equation}
We also know from \cite[Congruence (9)]{keith2022parity} that  
\begin{equation}\label{relfork=21(2)}
\frac{f_7}{f_1} \equiv f_1^6 + q f_1^2 f_7^4 +  q^2 \frac{f_7^8}{f_1^2}.
\end{equation}
From (\ref{relfork=21(1)}) and (\ref{relfork=21(2)}) we see that 
\[ 
\sum_{n = 0}^\infty \overline{C}_{84, 21} (n) q^n \equiv  \frac{f_{21}^3}{f_1} = \frac{f_{21}^3}{f_7} \frac{f_7}{f_1} 
\equiv 
\left( f_7^8 + q^7 \frac{f_{21}^{12}}{f_7^4} \right) 
\left(  f_1^6 + q f_1^2 f_7^4 +  q^2 \frac{f_7^8}{f_1^2}  \right). 
\] 
Extract the terms with even powers of $q$, then replace $q$ by $q^{1/2}$, and then multiply by $q$ to get
\begin{align*}
\sum_{n = 0}^\infty  \overline{C}_{84, 21} (2n) q^{n+1} &\equiv 
q f_1^3 f_7^4 + q^2\frac{f_7^{8}}{f_1} +  q^5 f_1 f_{21}^{6}  
\\
&= 
f_1 \left( q f_1^2 f_7^4 + q^2\frac{f_7^{8}}{f_1^2} \right) +  q^5 f_1 f_{21}^{6} .
\end{align*}
Using (\ref{relfork=21(2)}) we obtain
\begin{align*}
\sum_{n = 0}^\infty  \overline{C}_{84, 21} (2n) q^{n+1} 
&\equiv f_7 + f_1^7 +  q^5 f_1 f_{21}^{6} \\
&\equiv f_7 + f_1^3 f_4 +  q^5 f_1 f_{42}^{3}. 
\end{align*}
Each of the summands is lacunary modulo $2$ by congruences (\ref{euler}) to (\ref{jacobi-t}) and Landau's theorem. Therefore, the sequence $\overline{C}_{84, 21} (2n)$ is also lacunary modulo $2$. 
\end{proof}

\subsection{Case $k=23$}

\begin{proof}
We have
\[ \frac{f_{23}}{f_1} = \sum_{n = 0}^\infty b_{23}(n) q^n\]
where $b_{23}(n)$ is the number of 23-regular partitions of $n$. Congruence (4.4) in \cite{baruah2015parity} tells us that 
\[
\sum_{n = 0}^\infty b_{23}(2n+1) q^n \equiv f_{23}  + q f_1 f_{46}.
\]
By replacing $q$ with $q^2$, and multiplying the resulting congruence by $q$, we obtain
\[
\sum_{n = 0}^\infty b_{23}(2n+1) q^{2n+1} \equiv q f^2_{23}  + q^3 f_1^2 f_{46}^2.
\]
From this we find that the part of $f_{23}^3/f_1$ whose terms are those with odd powers of $q$ is given by
\[
\sum_{n = 0}^\infty \overline{C}_{92 , 23}(2n+1) q^{2n+1} \equiv 
q f_{23}^4 + q^3 f_1^2 f_{46}^2 f_{23}^2.
\]
Dividing by $q$ and replacing $q$ by $q^{1/2}$,
\[ 
\sum_{n = 0}^\infty \overline{C}_{92 , 23}(2n+1) q^{n} 
\equiv 
f_{23}^2 + q f_1 f_{46} f_{23} \equiv f_{23}f_{23} + q f_1 f_{23}^3.
\]
Using congruences (\ref{euler}), (\ref{euler-t}) and (\ref{jacobi-t}) and Landau's theorem, we find that each summand is lacunary modulo $2$. Thus, the sequence $\overline{C}_{92 , 23}(2n+1)$ is also lacunary modulo $2$.
\end{proof}

\printbibliography

\footnotesize
\begin{flushleft}
  \textsc{Institute of Mathematics, University of the Philippines Diliman, Quezon City, Philippines} \\
  \textit{Email:} \texttt{vmaricheta@math.upd.edu.ph}, \texttt{jdimabayao@math.upd.edu.ph}, \texttt{hjshi@math.upd.edu.ph}
\end{flushleft}

\end{document}